\documentclass[10pt]{amsart}

\usepackage{amsmath,amsthm,amssymb}

\newtheorem{thm}{Theorem}
\newtheorem{lem}[thm]{Lemma}
\newtheorem{prop}[thm]{Proposition}
\newtheorem{cor}[thm]{Corollary}

\theoremstyle{remark}
\newtheorem*{remark}{Remark}

\newcommand{\Z}{{\mathbb Z}}

\newcommand{\R}{{\mathbb R}}

\newcommand{\Prob}{{\mathbb{P}}} 
\newcommand{\<}{\left\langle}
\renewcommand{\>}{\right\rangle}

\newcommand{\intinf}{\int_{-\infty}^\infty}
\renewcommand{\d}{{\mathrm{d}}} 
\renewcommand{\i}{{\mathrm{i}}} 
\newcommand{\I}{1\!\!1} 

\renewcommand{\Im}{{\mathfrak{Im}}}
\renewcommand{\O}{{\mathcal{O}}} 
\renewcommand{\^}{\widehat} 

\newcommand{\smoothN}{\widetilde N}
\newcommand{\smoothS}{\widetilde S}
\newcommand{\Nf}{\smoothN_M} 
\newcommand{\Sf}{\smoothS_{M,L}} 

\newcommand{\mmt}{\mathcal{M}} 
\newcommand{\meas}{\operatorname{meas}}
\newcommand{\sqfree}{\ \square\text{free}}
\renewcommand{\th}{\textsuperscript{th}} 

\begin{document}
\title{On the distribution of lattice points in  thin annuli}
\author{C.P. Hughes}
\author{Z. Rudnick}
\address{Raymond and Beverly Sackler School of Mathematical
Sciences\\ Tel Aviv University\\ Tel Aviv 69978\\ Israel.}
\curraddr{American Institute of Mathematics\\ 360 Portage Avenue\\
Palo Alto, CA 94306-2244\\ U.S.A.} \email{hughes@aimath.org}
\address{Raymond and Beverly Sackler School of Mathematical
Sciences\\
Tel Aviv University\\ Tel Aviv 69978\\ Israel}
\email{rudnick@post.tau.ac.il}

\date{12 September 2003}

\begin{abstract}
We show that the number of lattice points lying in a thin annulus
has a Gaussian value distribution if the width  of the annulus
tends to zero sufficiently slowly as we increase the inner radius.
\end{abstract}

\maketitle

\section{Introduction}

Let $N(t)$ be the number of integer lattice points in a disk of
radius $t$ centered at the origin. Thus $N(t) = \sum_{n\leq t^2}
r(n)$ where $r(n)$ is the number of ways of writing $n=x^2+y^2$ as
a sum of two squares. As is well known, $N(t)$ is asymptotic to
the area $\pi t^2$ of the disk. Much effort has gone into
understanding the growth of the remainder term. Heath-Brown
\cite{H-B} considered the {\em distribution} of the normalized
remainder term $(N(t)-\pi t^2)/\sqrt{t}$, and proved that it has a
limiting value distribution in the sense that there exists a
probability distribution function $\nu$ such that for any interval
$\mathcal{A}$,
\begin{equation*}
\frac{1}{T} \meas\left\{ t \in [T,2T] : \frac{N(t)-\pi
t^2}{\sqrt{t}} \in \mathcal{A}\right\} \to \int_\mathcal{A}
\nu(x)\;\d x
\end{equation*}
where the measure is the ordinary Lebesgue measure. It is known
that $\nu(x)$ is not the Gaussian measure, for instance the tails
have been shown to decay roughly like $\exp(-x^4)$,
\cite{BCDL,HB2}.

Bleher, Dyson  and Lebowitz \cite{BDL, Ble-Leb, Ble} investigated
the distribution of a  similarly scaled remainder term of the
number $N(t,\rho):=N(t+\rho)-N(t)$ of lattice points in an annulus
of inner radius $t$ and width $\rho(t)$ depending on $t$. The
``expected'' number  of points is the area $\pi(2t\rho+\rho^2)$ of
the annulus. Define a normalized remainder term by
\begin{equation*}
S(t,\rho) := \frac{N(t+\rho)-N(t) - \pi(2t\rho+\rho^2)}{\sqrt{t}}
\;.
\end{equation*}

The picture that emerges is that there is a number of
distinct regimes:
\begin{itemize}
\item The ``global'', or ``macroscopic'', regime
$\rho(t)\to\infty$ (but $\rho=o(t)$), in which case Bleher and
Lebowitz \cite{Ble-Leb} show that $S(t,\rho)$ has a limiting
distribution with tails which decay roughly as $\exp(-x^4)$. In
fact the distribution is that of the difference of two i.i.d.
random variables whose distribution is the limiting distribution
of $(N(t)-\pi t^2)/\sqrt{t}$.

\item The intermediate, or ``mesoscopic'', regime $\rho \to 0$ but
$\rho t\to\infty$. The variance of $S(t,\rho)$ is given by
\cite{BL2}
$$
\frac 1T \int_T^{2T} |S(t,\rho)|^2 \d t \sim  \sigma^2:=16
\rho\log \frac 1\rho
$$
and  Bleher and Lebowitz~\cite{Ble-Leb} conjectured that
$S(t,\rho)/ \sigma$ has a standard Gaussian distribution.

\item The ``saturation'' regime $0<\rho(t)<\infty$ is fixed as
$t\to\infty$, where it has been shown \cite{Ble-Leb} that
$S(t,\rho)$ has a distribution with rapidly decaying tails. As
$\rho\to \infty$ the distribution converges to that found in the
macroscopic regime, and as $\rho\to0$ it converges to the
conjectured mesoscopic distribution. 

\item The local regime, $\rho \approx 1/t$. If  the annulus were
centered at a generic point rather than at a lattice point, or if
we consider ``generic'' lattices instead of the integer lattice
$\Z^2$, then it is consistent with conjectures of Berry and Tabor
\cite{BT} that the statistics are Poissonian; see \cite{Sarnak,
EMM, Marklof} for some progress on this, as well as \cite{Sinai,
Major, Minami}.
\end{itemize}

In this paper we prove part of the Gaussian distribution
conjecture of Bleher and Lebowitz. We will show that $S(t,\rho)$
has a Gaussian distribution when $\rho$ shrinks to zero
sufficiently slowly:

\begin{thm}\label{thm:Bleher conjecture}
If $\rho \to 0$ but $\rho\gg T^{-\delta}$ for all $\delta>0$, then
for any  interval $\mathcal{A}$
\begin{equation*}
\lim_{T\to\infty} \frac{1}{T}\meas\left\{t\in[T,2T]\ :\
\frac{S(t,\rho)}{\sigma} \in \mathcal{A}\right\} =
\frac{1}{\sqrt{2\pi}} \int_\mathcal{A} e^{-x^2/2}\;\d x
\end{equation*}
where $\sigma^2 = 16\rho \log \frac 1 {\rho}$.
\end{thm}

The structure of the argument is as follows: We replace the sharp
counting function $N(t)$ by a smooth counting function $\Nf(t)$
whose smoothness parameter $M=M(T)$ depends on $T$ (note that
though $t$ and $T$ are formally independent, we always think of
$t$ as being around $T$). Since we are only interested in $\rho\to
0$ we will set $\rho=1/L$ where $L=L(T)$ tends to infinity with
$T$, and we define the corresponding normalized remainder term to
be
\begin{equation*}
\Sf(t) := \frac{\Nf(t+1/L) - \Nf(t) - 2\pi/L - \pi/L^2}{\sqrt{t}}
\end{equation*}
We compute the moments of $\Sf(t)$ when $t$ is chosen at random
with respect to a smooth measure. We show in Section~\ref{sec: Nf}
that the $m$\textsuperscript{th} moment of $\Sf/\sigma$ converges
to that of a standard normal random variable provided $L \ll
T^{\nu(m)}$, with $0<\nu(m)<1/(2^{m-1}-1)$. Thus $\Sf$ has a
normal distribution if $L\to\infty$ but $L \ll T^{\delta}$ for all
$\delta>0$. In Section~\ref{sec: unsmoothing} we show that the
variance of the difference $(S(t,1/L)-\Sf(t))/\sigma$ goes to
zero, and hence $S(t,\rho)/\sigma$ has a normal distribution with
respect to the smooth measure. Finally we use an approximation
argument to pass from smooth measures to the Lebesgue measure used
in Theorem~\ref{thm:Bleher conjecture}.

\section{Smoothing}
To obtain Theorem~\ref{thm:Bleher conjecture} we will replace
sharp cutoffs by smooth ones. First, we will replace Lebesgue
measure with a smooth average of $t$ around $T$, that is we pick
$t$ at random by taking a smooth function $\omega\geq 0$, of total
mass unity, such that both $\omega$ and its Fourier transform
$\^\omega $ are rapidly decaying, in the sense that for any $A>2$,
\begin{equation}\label{eq:omega decay}
\omega(t) \ll \frac{1}{(1+|t|)^A}, \qquad \^\omega(t) \ll
\frac{1}{(1+|t|)^A}
\end{equation}
for all $t$. (In fact we also choose $\omega$ to be supported on
the positive reals as this makes the analysis simpler).

Define the averaging operator
\begin{equation}\label{eq:average}
\< f \> = \frac{1}{T} \intinf f(t)
\omega\left(\frac{t}{T}\right)\;\d t
\end{equation}
(this is the expected value of $f$ with respect to this measure),
and let $\Prob_{\omega,T}$ be the associated probability measure:
$$\Prob_{\omega,T}(f\in \mathcal A) =
\frac{1}{T} \intinf \I_{\mathcal A}( f(t))
\omega\left(\frac{t}{T}\right)\;\d t
$$
(Throughout the paper we will extend $N(t)$, $S(t,\rho)$ and
similar functions, initially defined for $t>0$, to the whole real
line. Since $\omega(t)=0$ for $t\leq 0$ we are free to choose
whichever extension makes the analysis most simple).

We will also smooth the edges of the circle, and show that this
modified counting function has a Gaussian distribution: Let $\chi$
be the indicator function of the unit disc, and $\psi$ a smooth,
even, function on the real line, of total mass unity, whose
Fourier transform $\^\psi$ is smooth and has compact support.
Define a rotationally symmetric function $\Psi$ on $\R^2$ by
setting $\^\Psi(\vec y)=\^\psi(|\vec y|)$ where $|\vec y|$ denotes
the standard Euclidean norm of $\vec y\in \R^2$, and where the
Fourier transform is

\begin{equation*}
\^ f(\vec{y}) = \int_{\R^2} f(\vec{x})
e^{-2\pi\i\<\vec{x},\vec{y}\>}\;\d \vec{x}
\end{equation*}
with $\<\vec{x},\vec{y}\>$ the usual Euclidean inner product. For
$\epsilon>0$ set
\begin{equation}\label{eq:defn Psi}
\Psi_\epsilon(\vec x) = \frac 1{\epsilon^2}\Psi(\frac {\vec
x}{\epsilon})
\end{equation}
Now set $\chi_\epsilon = \chi *\Psi_\epsilon$ to be the
convolution of $\chi$ and $\Psi_\epsilon$, which is a smoothed
indicator function of the unit disc with ``fuzziness'' of width
$\epsilon$, in the sense that $0\leq \chi_\epsilon \leq 1$, and if
$\psi$ (rather than its Fourier transform $\^\psi$) had compact
support then $\chi-\chi_\epsilon$ would be concentrated in the
shell $1-\epsilon<|\vec x|<1+\epsilon$. Due to the rapid decay of
tails, this is essentially still the case when $\psi$ is in the
Schwarz class, as it is for us.

Now take $\epsilon = 1/t\sqrt{M}$ where $M=M(T)$ depends on $T$
and tends to infinity with $T$, and  define a smooth counting
function, or smooth linear statistic, by
$$
\Nf(t) = \sum_{\vec n\in \Z^2} \chi_\epsilon(\frac {\vec n}t)
$$
This counts lattice points in a ``fuzzy circle'' of radius about
$t$, with fuzziness about $t\epsilon = 1/\sqrt{M}$.

The number of lattice points in a smooth annulus of inner radius
$t$ and width $\rho$ is therefore given by $\Nf(t+\rho)-\Nf(t)$.
Since we are interested in radii $t$ in an interval $[T,2T]$,  we
will in what follows freeze the width of the annulus to be
$\rho(T)$ as $t$ varies in $[T,2T]$ rather than allowing it to
vary  with $t$; this will simplify some of the calculations.
Furthermore, since from henceforth we are only concerned with
$\rho\to 0$, we will set $\rho=1/L$, and let $L(T)\to\infty$ as
$T\to \infty$.

Set
\begin{equation}\label{eq:defn Sf}
\Sf = \frac{ \Nf(t+1/L) - \Nf(t) - 2\pi t /L-\pi/L^2}{\sqrt{t}}
\end{equation}
The width of the smoothed sides of $\Nf$ is $\O(\epsilon t) =
\O(1/\sqrt{M})$. In order for $\Sf$ to approximate $S(t,1/L)$, it
must be that $1/L$ is much larger than the width of the sides, so
we insist that $ L/\sqrt{M} \to 0$.

We will show:
\begin{thm}\label{thm:distrn for Sf}
Suppose that  $M(T)$, $L(T)$ are increasing to infinity with $T$
such that $M= \O(T^\delta)$ for all $\delta>0$, and  $L/\sqrt{M}
\to 0$, then for any  interval $\mathcal{A}$,
\begin{equation*}
\lim_{T\to\infty} \Prob_{\omega,T} \left\{ \frac{\Sf}{\sigma} \in
\mathcal{A}\right\} = \frac{1}{\sqrt{2\pi}} \int_{\mathcal{A}}
e^{-x^2/2}\;\d x
\end{equation*}
where $\Sf$ is given by \eqref{eq:defn Sf} and
\begin{equation*}
\sigma^2 = \frac{16\log L}{L}
\end{equation*}
\end{thm}

\begin{remark}
The arguments given below for the proof of Theorem~\ref{thm:distrn
for Sf} will also prove a central limit theorem for smooth linear
statistics in higher dimensions: Defining $\chi_\epsilon=\chi \ast
\Psi_\epsilon$ where $\chi$ is the indicator function of the unit
ball and $\Psi_\epsilon$ is defined in analogy with \eqref{eq:defn
Psi}, we have a smooth counting function $\Nf(t) := \sum_{\vec
n\in \Z^d} \chi_\epsilon(\frac {\vec n}t) \;.$ where as before
$\epsilon=1/t\sqrt{M}$.

The asymptotic behaviour of $\Nf(t)$ is given by $c_d t^d$, with
$c_d$ the volume of the unit ball in $\R^d$.

It may then be shown that if $M=\O(T^\delta)$ for all $\delta>0$,
then the distribution of the normalized remainder term
$\smoothS_{M}(t)= \frac{\Nf(t) - c_d t^d}{t^{(d-1)/2}} $ when
averaged over $t$ around $T$ weakly converges to a Gaussian with
mean zero and variance
$$
\sigma^2 = \begin{cases}
\frac{2}{\pi^2} K_3 \log M & \text{ when } d=3\\
\frac{d-1}{\pi^2} K_d \int_0^\infty y^{d-4} \^\psi(y)^2\;\d y
M^{(d-3)/2} & \text{ when } d\geq 4
\end{cases}\label{eq:asym sigma high dims}
$$
where
\begin{equation*}
K_d = \frac{4^{d-1} \pi^{d-1/2}}{2^d - 1}
\frac{\Gamma(\tfrac12d-\tfrac12)}{\Gamma(d) \Gamma(\tfrac12d)}
\frac{\zeta(d-1)}{\zeta(d)} \;.
\end{equation*}
\end{remark}

\section{The distribution of $\Nf$}\label{sec: Nf}

\begin{lem}\label{lem:Nf Poisson}
As $t\to\infty$, we have
\begin{equation*}
\Nf(t) = \pi t^2 -\frac{\sqrt{t}}{\pi} \sum_{n=1}^\infty
\frac{r(n)}{n^{3/4}}  \cos(2\pi t \sqrt{n} +\tfrac14\pi)
\^\psi(\sqrt{\frac nM}) + \O\left(\frac{1}{\sqrt{t}}\right)
\end{equation*}
with the error term independent of $M$.
\end{lem}

\begin{proof}
By Poisson summation,
\begin{align}
\Nf(t) &:= \sum_{\vec{n} \in \Z^2} (\chi*\Psi_{\epsilon})\left(\frac{\vec{n}}{t}
\right) \nonumber\\
&= t^2\sum_{\vec{k}\in\Z^2} \^ \chi(t \vec{k}) \^\Psi_{\epsilon}(t
\vec{k}) \label{eq:Nf Fourier}
\end{align}

Changing into polar coordinates, and using the fact that $\chi$ is
rotationally symmetric, the 2-dimensional Fourier transform of
$\chi$ is
\begin{align*}
\^\chi(\vec{y}) &= \int_0^1 r \int_0^{2\pi} e^{-2\pi\i r |\vec{y}| \cos\theta}
\;\d\theta\;\d r \\
&= \frac{-\cos(2\pi |\vec{y}|+\tfrac14\pi)}{\pi |\vec{y}|^{3/2}} +
\O\left(\frac{1}{|\vec{y}|^{5/2}} \right)
\end{align*}
as $|\vec{y}|\to\infty$. By its definition in \eqref{eq:defn Psi},
$\^{\Psi}_\epsilon(\vec{y}) = \^\Psi(\epsilon
\vec{y})=\^\psi(\epsilon |\vec{y}|)$. Therefore, inserting this
into \eqref{eq:Nf Fourier}, treating the mean (when
$\vec{k}=\vec{0}$) separately, and setting $\epsilon=1/t\sqrt{M}$,
\begin{align*}
\Nf(t) &= \pi t^2 - \frac{\sqrt{t}}{\pi} \sum_{\vec{k} \neq
\vec{0}} \left\{ \frac{\cos(2\pi t
|\vec{k}|+\tfrac14\pi)}{|\vec{k}|^{3/2}} \^\psi(\epsilon t
|\vec{k}|) + \O\left(\frac{1}{t} \frac{\^\psi(\epsilon t
|\vec{k}|)}{|\vec{k}|^{5/2}}\right)\right\}\\
&= \pi t^2 -\frac{\sqrt{t}}{\pi} \sum_{n=1}^\infty
\frac{r(n)}{n^{3/4}}  \cos(2\pi t \sqrt{n} +\tfrac14\pi) \^\psi(
\sqrt{\frac nM}) + \O\left(\frac{1}{\sqrt{t}}\right)
\end{align*}
with the constant implicit in the error term independent of
$M(T)$.
\end{proof}
Note that the compact support of $\^\psi$ means that the sum
truncates at $n\approx M$. Thus we need $M\gg 1$ in order for
there to be any terms in the sum.

Now, since
\begin{equation*}
\Sf = \frac{\Nf(t+1/L) - \Nf(t) - \pi(2t/L+1/L^2)}{\sqrt{t}}
\end{equation*}
then for $t\geq 1$ and $L\geq 1$,
\begin{align}
\Sf &= \frac{1}{\pi} \sum_{n=1}^\infty \frac{r(n)}{n^{3/4}}
\Bigl[\cos(2\pi t \sqrt{n} +\tfrac\pi 4) - \cos(2\pi
(t+\tfrac{1}{L}) \sqrt{n} +\tfrac\pi 4)\Bigr]
\^\psi(\sqrt{\frac{n}{M}}) + \O\left(\frac{1}{t}\right) \nonumber\\
&=\frac{2}{\pi} \sum_{n=1}^\infty \frac{ r(n)}{ n^{3/4}}
\sin\left(\frac{\pi \sqrt{n}}{L}\right) \sin\left(2\pi
(t+\tfrac{1}{2L})\sqrt{n}+\tfrac\pi 4\right)
\^\psi(\sqrt{\frac{n}{M}})  + \O\left(\frac{1}{t}\right)
\label{eq:smooth S infinite sum}
\end{align}

Note that we have three independent variables. The variable $t$,
which we always consider to be large, is the radius of the
annulus. This is the variable we average over. The width of the
annulus is $1/L$. Since we want a thin annulus, $L\to \infty$, and
Gaussian behaviour is not seen if this condition does not hold.
The annulus does not have sharp sides, but smoothed edges, and the
third independent variable is $M$; the larger $M$ is, the sharper
the annulus' sides (in the sense that it better approximates the
indicator function). We must have $L/\sqrt{M}\to0$ in order for
the annulus to have some width, and not be ``just sides''. That
is, the annulus shouldn't be too smooth.

\begin{proof}[Proof of Theorem~\ref{thm:distrn for Sf}]
First we show the mean is $\O(1/T)$. Since $\omega(t)$ is real,
\begin{equation*}
\< \sin\left(2\pi (t+\tfrac{1}{2L})\sqrt{n}+\tfrac14\pi\right) \>
= \Im\left\{ \^\omega(-T\sqrt{n})
e^{\i\pi(\frac{\sqrt{n}}{L}+\frac14)} \right\} \ll \frac{1}{T^A
n^{A/2}}
\end{equation*}
for any $A>2$, where we have used the rapid decay of $\^\omega$.
Thus
\begin{align*}
\<\Sf\> &\ll \sum_{n=1}^\infty \frac{ r(n)}{ n^{3/4}} \frac{1}{T^A
n^{A/2}} + \O(\frac{1}{T})\\
&=\O(\frac{1}{T})
\end{align*}

Setting
\begin{equation}\label{eq:defn M_k}
\mmt_m := \<\left(\frac{2}{\pi} \sum_{n=1}^\infty \frac{ r(n)}{
n^{3/4}} \sin\left(\frac{\pi \sqrt{n}}{L}\right) \sin\left(2\pi
(t+\tfrac{1}{2L})\sqrt{n}+\tfrac14\pi\right)
\^\psi(\sqrt{\frac{n}{M}})\right)^m\>
\end{equation}
then from \eqref{eq:smooth S infinite sum}, the Cauchy-Schwartz
inequality implies that the $m\th$ moment of $\Sf$ is
\begin{multline}\label{eq:mth moment Sf}
\< (\Sf)^m\> =\\
\<\left\{\frac{2}{\pi} \sum_{n=1}^\infty \frac{ r(n)}{ n^{3/4}}
\sin\left(\frac{\pi \sqrt{n}}{L}\right) \sin\left(2\pi
(t+\tfrac{1}{2L})\sqrt{n}+\tfrac14\pi\right)
\^\psi(\sqrt{\frac{n}{M}})+\O\left(\frac{1}{T}\right)\right\}^m\> \\
= \mmt_m + \O\left(\sum_{j=1}^m \binom{m}{j}
\frac{\sqrt{\mmt_{2m-2j}} }{T^j}\right)
\end{multline}

The conditions of Theorem~\ref{thm:distrn for Sf} are that
$M=\O(T^\delta)$ for all $\delta>0$, and that $L\to\infty$ in such
a way that $L/\sqrt{M}\to 0$. In which case,
Proposition~\ref{lem:variance} allows us to deduce that $\sigma^2
:= \mmt_2 \sim \frac{16\log L}{L}$, and Proposition~\ref{lem:mmts
Gaussian} shows that for all $m>2$,
\begin{equation*}
\frac{\mmt_m}{\sigma^m} =
\begin{cases}
\frac{m!}{2^{m/2} (m/2)!} + \O\left(\frac{1}{L^{1-\delta'}}\right) & \text{ if
$m$ is even}\\
\O\left(\frac{1}{L^{1-\delta'}}\right) & \text{ if $m$ is odd}
\end{cases}
\end{equation*}

These are the moments of the standard normal distribution, and
inserting these into \eqref{eq:mth moment Sf}, we see this is
sufficient to prove that the distribution of $\Sf/\sigma$ weakly
converges as $T\to\infty$ to a Gaussian with mean zero and
variance $1$.
\end{proof}

\subsection{The variance}\label{sect:variance}

\begin{prop}\label{lem:variance}
If $M = \O(T^{2(1-\delta)})$ for fixed $\delta>0$, then the
variance of $\Sf$ is asymptotic to
\begin{equation}\label{eq:sigma}
\sigma^2 := \frac{2}{\pi^2} \sum_{n=1}^\infty
\frac{r(n)^2}{n^{3/2}} \sin^2\left(\frac{\pi\sqrt{n}}{L}\right)
\^\psi^2\left(\sqrt{\frac{n}{M}}\right)
\end{equation}
If $L\to \infty$ but $L/\sqrt{M} \to 0$, then
\begin{equation}\label{eq:sigma_asympt}
\sigma^2 \sim \frac{16\log L}{L}
\end{equation}
\end{prop}

\begin{proof}
Expanding out \eqref{eq:defn M_k} we have
\begin{multline*}
\mmt_2 = \frac{4}{\pi^2} \sum_{m,n} \frac{r(m) r(n) \sin(\pi
\sqrt{m}/L) \sin(\pi \sqrt{n}/L) \^\psi(\sqrt{\frac{m}{M}})
\^\psi(\sqrt{\frac{n}{M}})}{(mn)^{3/4}} \times\\
\times
\<\sin\left(2\pi(t+\tfrac{1}{2L})\sqrt{m}+\tfrac14\pi\right)\sin\left(2\pi
(t+\tfrac{1}{2L})\sqrt{n}+\tfrac14\pi\right)\>
\end{multline*}

Now the average on the bottom line is
\begin{multline*}
\frac{1}{4} \Bigl[ \^\omega\left(T(\sqrt{m}-\sqrt{n})\right)
e^{\i\pi\frac1L(\sqrt{n}-\sqrt{m})} +
\^\omega\left(T(\sqrt{n}-\sqrt{m})\right)
e^{\i\pi\frac1L(\sqrt{m}-\sqrt{n})}\\
-
\^\omega\left(T(\sqrt{m}+\sqrt{n})\right)e^{-\i\pi(\frac12+\frac1L(\sqrt{m}
+\sqrt{n}))}-
\^\omega\left(-T(\sqrt{m}+\sqrt{n})\right)e^{\i\pi(\frac12+\frac1L(\sqrt{m}
+\sqrt{n}))}
\Bigr]
\end{multline*}
The support condition on $\^\psi$ means that $m,n$ are both
constrained to be $\O(M)$, and so either $m=n$ or else
$|\sqrt{m}-\sqrt{n}| \gg \frac{1}{\sqrt{M}}$. Using the bound
$\^\omega(t) \ll (1+|t|)^{-A}$ for all $A>0$, the off-diagonal
terms contribute at most
\begin{equation*}
\sum_{1\leq n\neq m\leq M} \left(\frac{\sqrt{M}}{T}\right)^A \ll
\frac{M^{A/2+2}}{T^A} \ll T^{4-\delta A}
\end{equation*}
using the assumption $M=\O(T^{2(1-\delta)})$. Therefore, for any
$B>0$,
\begin{equation*}
\mmt_2 = \frac{2}{\pi^2} \sum_{n=1}^\infty \frac{r(n)^2}{n^{3/2}}
\sin^2\left(\frac{\pi \sqrt{n}}{L}\right)
\^\psi^2\left(\sqrt{\frac{n}{M}}\right) + \O(T^{-B})
\end{equation*}
Define $\sigma^2$ to be the infinite sum above. Since $r(n)\ll
n^\eta$ for all $\eta>0$, $\sigma^2$ is bounded for all $L$. To
find the asymptotics as  $L\to\infty$, we use a formula of
Ramanujan \cite{Ramanujan}\footnote{We wish to thank Bob Vaughan
for pointing this reference out to us.}:
\begin{equation*}
\sum_{n \leq X} r(n)^2 = 4 X\log X + O(X)\;.
\end{equation*}

We then have
\begin{align*}
\sigma^2 &:= \frac{2}{\pi^2} \sum_{n=1}^\infty
\frac{r(n)^2}{n^{3/2}} \sin^2\left(\frac{\pi \sqrt{n}}{L}\right)
\^\psi^2\left(\sqrt{\frac{n}{M}}\right)\\
&\sim \frac{8}{\pi^2} \int_1^\infty \frac{\log x}{x^{3/2}}
\sin^2(\frac{\pi  \sqrt{x}}{L}) \^\psi^2(\sqrt{\frac{x}{M}})\;\d x\\
&= \frac{32}{L \pi^2} \int_{1/L}^\infty \log (y L)
\frac{\sin^2(\pi y)}{y^2} \^\psi^2\left(\frac{y
L}{\sqrt{M}}\right)\;\d y\\
&\sim \frac{\log L}{L} \frac{32}{\pi^2} \int_0^\infty
\frac{\sin^2(\pi y)}{y^2} \^\psi^2\left(\frac{y
L}{\sqrt{M}}\right)\;\d y
\end{align*}
on changing variables to $x=y^2 L^2$, and using the fact that we
assume that $L\to \infty$. Now using the additional restriction
(caused by the fuzziness of the annulus' sides) that $L/\sqrt{M}
\to 0$, we see that since $\^\psi(y L /\sqrt{M}) \sim 1$ for all
$y =o(\sqrt{M}/L)$, the integral can be evaluated asymptotically
to equal $\pi^2/2$, and so
\begin{equation*}
\sigma^2 \sim \frac{16\log(L)}{L}
\end{equation*}
Since $L=o(T^{1-\delta})$, the error terms in \eqref{eq:mth moment
Sf} are all smaller than $\sigma^2$, and so the variance of $\Sf$
is asymptotic to $\sigma^2$ as $T\to\infty$.
\end{proof}

The constraints on $M$, that $M=\O(T^{2-2\delta})$ but
$L/\sqrt{M}\to 0$, illustrate the role of smoothing. The first
constraint, that $M$ is not too big, comes from requiring that the
annulus is sufficiently smooth to easily handle the averages (to
enable us to reduce to the diagonal). The second constraint, that
$M$ is not too small, is to ensure the function is not too smooth,
so that the width of the edges is greater than the size of the
annulus. (That $L\to \infty$ forces $M$ to go to infinity. If it
did not, the function would be so smooth as to have no
fluctuations!)

\subsection{The higher moments}

\begin{prop}\label{lem:mmts Gaussian}
For fixed $\delta>0$, if
$M=\O\left(T^{2(1-\delta)/(2^{m-1}-1)}\right)$, and if $L\to
\infty$ such that $L/\sqrt{M}\to 0$, then for arbitrary
$\delta'>0$,
\begin{equation*}
\frac{\mmt_m}{\sigma^m} =
\begin{cases}
\frac{m!}{2^{m/2} (m/2)!} + \O\left(\frac{1}{L^{1-\delta'}}\right) & \text{ if
$m$ is even}\\
\O\left(\frac{1}{L^{1-\delta'}}\right) & \text{ if $m$ is odd}
\end{cases}
\end{equation*}
where $\mmt_m$ is given in \eqref{eq:defn M_k} and $\sigma^2$ is
given in \eqref{eq:sigma}.
\end{prop}

We will need to give lower bounds for alternating sums $\sum\pm
\sqrt{n_j}$. To do so, we use  the following lemma, a form of
Liouville's theorem, (cf \cite{H-B}):
\begin{lem}\label{lem:HB}
For $j=1,\dots,m$, let $n_j\leq M$ be positive integers. Then
\begin{itemize}
\item either $\sum \epsilon_j \sqrt{n_j} = 0$ for some $\epsilon_j=\pm1$,
\item or for all $\epsilon_j=\pm1$,
\begin{equation*}\left| \sum_{j=1}^m \epsilon_j \sqrt{n_j} \right|
\geq \frac{1}{(m\sqrt{M})^{2^{m-1}-1}}\end{equation*}
\end{itemize}
\end{lem}

\begin{proof}
Assume that $\sum \epsilon_j \sqrt{n_j} \neq 0$ for all choices of
$\epsilon_j=\pm1$. Then
\begin{equation*}
P := \prod_{\epsilon_j = \pm 1} \left( \sum_{j=1}^m \epsilon_j
\sqrt{n_j}\right)
\end{equation*}
is non-zero. By Galois theory, since $\sum\epsilon_j\sqrt{n_j}$ is
an algebraic number, and $P$ is the product over all possible
symmetries, $P$ is an integer. Since we assumed that no term in
$P$ vanishes, $|P|\geq 1$. Since both $\sum \epsilon_j \sqrt{n_j}$
and $-\sum \epsilon_j \sqrt{n_j}$ are terms in $P$, if
\begin{equation*}
Q := \prod_{\substack{\epsilon_j = \pm 1\\j=2,3,\dots,m}} \left(
\sqrt{n_1}+\sum_{j=2}^m \epsilon_j \sqrt{n_j}\right)
\end{equation*}
then $P=(-1)^{2^{m-1}} Q^2$, and so $|Q| = \sqrt{|P|} \geq 1$.

By assumption $n_j \leq M$ for all $j$, and so, independent of the
$\epsilon_j$,
\begin{equation*}
\left|\sqrt{n_1}+\sum_{j=2}^m \epsilon_j \sqrt{n_j}\right| \leq m
\sqrt{M}
\end{equation*}
and so for any $\eta_j = \pm 1$,
\begin{equation*}
\left|\sqrt{n_1} + \sum_{j=2}^m \eta_j \sqrt{n_j}\right| =
\frac{|Q|}{\prod^* \left|\sqrt{n_1} + \sum_{j=2}^m \epsilon_j
\sqrt{n_j}\right|} \geq \frac{1}{(m \sqrt{M})^{2^{m-1}-1}}
\end{equation*}
where $\prod^*$ denotes the product over all $\epsilon_j$ distinct
from $\eta_j$, there being $2^{m-1}-1$ terms in such a product.
\end{proof}

{}From this, it is simple to derive the following lemma.
\begin{lem}\label{lem:sqrt n_j zero or large}
For $j=1,\dots,m$ let $n_j\leq M$ be positive integers, and let
$\epsilon_j = \pm 1$ be such that
\begin{equation*}
\sum_{j=1}^m \epsilon_j \sqrt{n_j} \neq 0
\end{equation*}
Then
\begin{equation*}
\left| \sum_{j=1}^m \epsilon_j \sqrt{n_j} \right| \geq
\frac{1}{(m\sqrt{M})^{2^{m-1}-1}}
\end{equation*}
\end{lem}

\begin{proof}
Either $\sum \eta_j \sqrt{n_j} \neq 0$ for any choice of
$\eta_j=\pm1$, and we are done by Lemma~\ref{lem:HB}, or else
there exists a (strict) subset $S \subsetneq \{1,\dots,m\}$ such
that
\begin{equation*}
\sum_{j \in S} \epsilon_j \sqrt{n_j} - \sum_{j \not\in S}
\epsilon_j \sqrt{n_j} = 0
\end{equation*}
so that
\begin{equation*}
\left|\sum_{j=1}^m \epsilon_j \sqrt{n_j}\right| = 2\left|\sum_{j
\in S} \epsilon_j \sqrt{n_j}\right|
\end{equation*}
Note that, by assumption, $\sum_{j\in S} \epsilon_j \sqrt{n_j}
\neq 0$ and, if $m'$ denotes the number of terms in the sum, then
$1\leq m' < m$. Now repeat the argument: Either $\sum_{j\in S}
\eta_j \sqrt{n_j} \neq 0$ for any choice of $\eta_j=\pm1$, in
which case Lemma~\ref{lem:HB} gives that
\begin{equation*}
\left| \sum_{j \in S} \epsilon_j \sqrt{n_j} \right| \geq
\frac{1}{(m'\sqrt{M})^{2^{m'-1}-1}} >
\frac{1}{(m\sqrt{M})^{2^{m-1}-1}}
\end{equation*}
or else one can find a further subdivide the set $S$ as before.
Since the number of terms in the sum is a positive integer and
reduces upon each subdivision, this process terminates.
\end{proof}

\begin{proof}[Proof of Proposition~\ref{lem:mmts Gaussian}.]
Expanding \eqref{eq:defn M_k} out,
\begin{multline*}
\mmt_m = \frac{2^m}{\pi^m} \sum_{n_1,\dots,n_m \geq 1}
\prod_{j=1}^m \frac{ r(n_j)}{ n_j^{3/4}} \sin\left(\frac{\pi
\sqrt{n_j}}{L}\right) \^\psi(\sqrt{\frac{n_j}{M}})\\
\times \< \prod_{j=1}^m\sin\left(2\pi
(t+\frac{1}{2L})\sqrt{n_j}+\tfrac14\pi\right)\>
\end{multline*}

Now,
\begin{multline*}
\< \prod_{j=1}^m \sin\left(2\pi(t+\tfrac{1}{2L})\sqrt{n_j}+\tfrac14\pi\right)
\> \\
= \< \prod_{j=1}^m \frac{1}{2\i}
\left[e^{2\pi\i((t+1/2L)\sqrt{n_j}+1/8)} -
e^{-2\pi\i((t+1/2L)\sqrt{n_j}+1/8)} \right]\>\\
=\sum_{\epsilon_j = \pm1} \frac{1}{2^m \i^m} \intinf \prod_{j=1}^m
\epsilon_j \exp\left( \epsilon_j
2\pi\i\left((t+\frac{1}{2L})\sqrt{n_j}+\tfrac18\right)\right)
\frac{1}{T} \omega\left(\frac{t}{T}\right)\;\d t\\
= \sum_{\epsilon_j = \pm1} \frac{\prod \epsilon_j}{2^m \i^m}
\^\omega\left(-T \sum_{j=1}^m \epsilon_j\sqrt{n_j}\right)
e^{\sum_{j=1}^m \epsilon_j \pi
\i\left(\frac{1}{L}\sqrt{n_j}+1/4\right)}
\end{multline*}

By the compact support condition of $\^\psi$, we may always assume
that $n_j =\O(M)$. By Lemma~\ref{lem:sqrt n_j zero or large} and
the fact that $\^\omega$ decays faster than any polynomial power,
the off-diagonal terms (those terms with $\sum_{j=1}^m \epsilon_j
\sqrt{n_j} \neq 0$) contribute at most
\begin{equation*}
\sum_{1\leq n_1,\dots,n_m\leq M}
\left(\frac{(\sqrt{M})^{2^{m-1}-1}}{T}\right)^A \ll
\frac{M^{(2^{m-1}-1)A/2+m}}{T^A} \ll T^{-\delta A +
2m/(2^{m-1}-1)}
\end{equation*}
which is vanishingly small, since $A$ can be arbitrarily large.
Thus the only contributing terms are those with $\sum_{j=1}^m
\epsilon_j \sqrt{n_j} = 0$, and using the fact that
$\^\omega(0)=1$, we therefore have for any $B>0$,
\begin{equation*}
\mmt_m =  \sum_{n_1,\dots,n_m} \sum_{\substack{\epsilon_j = \pm
1\\\sum \epsilon_j \sqrt{n_j} = 0}} \prod_{j=1}^m
\frac{-\i\epsilon_j r(n_j)}{\pi n_j^{3/4}} \sin\left(\frac{\pi
\sqrt{n_j}}{L}\right) \^\psi(\sqrt{\frac{n_j}{M}}) e^{\i
\pi\epsilon_j/4} + \O(T^{-B})
\end{equation*}

In order to estimate the size of $\mmt_m/\sigma^m$ when $L\to
\infty$, we need to use a lemma of  Besicovitch~\cite{Bes}.
\begin{lem}
If $q_j$, (for $j=1,\dots,m$), are distinct squarefree positive
integers, then $\sqrt{q_1},\dots,\sqrt{q_m}$ are linearly
independent over the rationals.
\end{lem}

Therefore, if $\sum_{j=1}^m \epsilon_j \sqrt{n_j} = 0$ with $n_j
\geq 1$, then there must exists a division of $\{1,\dots,m\}$ into
$\{S_i\}$ such that
\begin{equation*}
\{1,\dots,m\} = \coprod_{i=1}^\ell S_i
\end{equation*}
where $\sum_{i=1}^\ell |S_i|=m$ such that for $i=1,2,\dots,\ell$,
for all $j\in S_i$, $n_j = q_i f_j^2$ with the $q_i$ being
distinct square-free integers, and with the $f_j$ satisfying
\begin{equation*}
\sum_{j \in S_i} \epsilon_j f_j = 0
\end{equation*}

Summing over all possible divisions we see that
\begin{multline*}
\frac{\mmt_m}{\sigma^m} = \sum_{\ell=1}^m \sum_{\{1,\dots,m\} =
\coprod_{i=1}^\ell S_i} \Biggl(\frac{1}{\sigma^{|S_1|}}\sum_{q_1
\sqfree} D_{q_1}(S_1)\Biggr) \\
\times\Biggl(\frac{1}{\sigma^{|S_2|}}\sum_{\substack{q_2
\sqfree\\q_2\neq q_1}} D_{q_2}(S_2)\Biggr) \cdots
\Biggl(\frac{1}{\sigma^{|S_\ell|}}\sum_{\substack{q_\ell
\sqfree\\q_\ell\neq q_1,\dots,q_{\ell-1}}}
D_{q_\ell}(S_\ell)\Biggr)
\end{multline*}
where
\begin{multline}\label{eq:defn D_q(S)}
D_q(S) := \\
\frac{1}{q^{3|S|/4}} \sum_{\substack{f_j \geq 1\\\epsilon_j = \pm
1\\\sum_{j\in S} \epsilon_j f_j = 0}} \prod_{j\in S}
\frac{-\i\epsilon_j e^{\i\pi\epsilon_j /4} r(q f_j^2)}{\pi
f_j^{3/2}} \sin\left(\pi \tfrac{1}{L} f_j \sqrt{q}\right)
\^\psi\left(f_j \sqrt{\frac{q}{M}}\right)
\end{multline}

We will show in Lemma~\ref{lem:size of fundamental block} that if
$L \to \infty$ such that $L/\sqrt{M}\to 0$ then for all
$\delta'>0$,
\begin{equation*}
\frac{1}{\sigma^{|S|}}\sum_{q \sqfree} D_q(S) =
\begin{cases}
0 & \text{ if } |S|=1\\
1 & \text{ if } |S|=2\\
\O\left(\frac{1}{L^{1-\delta'}}\right) & \text{ otherwise}
\end{cases}
\end{equation*}

Therefore the only terms in $\mmt_m / \sigma^m$ which do not
vanish as $L \to \infty$ are those where $|S_i| = 2$ for all $i$.
If $m$ is odd, there are no such terms, and if $m=2k$ is even,
then the number of terms is equal to the number of ways of
partitioning $\{1,\dots,2k\}$ into $\coprod_{i=1}^k S_i$ with
$|S_i|=2$, which equals
\begin{equation*}
\frac{1}{k!} \binom{2k}{2} \binom{2k-2}{2} \dots \binom{2}{2} =
\frac{(2k)!}{k! 2^k}
\end{equation*}
This completes the proof of Proposition~\ref{lem:mmts Gaussian}
\end{proof}

\begin{lem} \label{lem:size of fundamental block}
If $L \to \infty$ is such that $L/\sqrt{M}\to0$ then
\begin{equation*}
\frac{\sum_{q \sqfree} D_q(S)}{\sigma^{|S|}} =
\begin{cases}
1 & \text{ if } |S|=2\\
\O\left(\frac{1}{L^{1-\delta}}\right) & \text{ otherwise}
\end{cases}
\end{equation*}
for all $\delta>0$, where $D_q(S)$ is defined in \eqref{eq:defn
D_q(S)}, and $\sigma^2$ is defined in \eqref{eq:sigma}.
\end{lem}

\begin{proof}
For convenience we assume, without loss of generality, that $S =
\{1,2,\dots,|S|\}$. Using $r(n)\ll n^{\delta}$ for all $\delta>0$,
and $\^\psi(x) \ll 1$, we can upper bound by
\begin{equation}\label{eq:upper bound D1}
\sum_{q \sqfree} D_q(S) \ll \sum_{q=1}^\infty \frac{q^{|S|
\delta}}{q^{3|S|/4}} Q(q)
\end{equation}
where
\begin{equation*}
Q(q) = \sum_{\epsilon_j=\pm1} \sum_{\substack{ f_j \geq 1 \\
\sum_{j=1}^{|S|} \epsilon_j f_j=0}} \prod_{j=1}^{|S|}
\frac{f_j^{\delta}}{f_j^{3/2}} \left|\sin\left(\pi f_j
\sqrt{q}/L\right)\right|
\end{equation*}

Note that $Q(q) \ll 1$ for all $q$. When $q \ll L^2$ a sharper
result can be deduced by a more careful treatment of $Q(q)$. In
order for $\sum_{j=1}^{|S|} \epsilon_j f_j = 0$, at least two of
the $\epsilon$ must have different signs, and so, with no loss of
generality, we put $\epsilon_{|S|}=-1$ and $\epsilon_{|S|-1}=+1$.
Hence
\begin{equation*}
f_{|S|} = f_{|S|-1} +  \sum_{j=1}^{|S|-2} \epsilon_j f_j
\end{equation*}
In order for both $f_{|S|} \geq 1$ and $f_{|S|-1}\geq 1$, it must
be that
\begin{equation*}
f_{|S|-1}\geq 1+\max\left\{0,-\sum_{j=1}^{|S|-2} \epsilon_j
f_j\right\}
\end{equation*}
Therefore,
\begin{multline*}
Q(q) = 2\sum_{\epsilon_1,\dots,\epsilon_{|S|-2}=\pm1}
\sum_{f_1,\dots,f_{|S|-2} \geq 1} \sum_{f_{|S|-1} \geq
1+\max\{0,-\sum_{j=1}^{|S|-2} \epsilon_j f_j\}} \\
\left(\prod_{j=1}^{|S|-1}  \frac{\left|\sin\left(\pi f_j
\sqrt{q}/L\right)\right|}{f_j^{3/2-\delta}}\right)  \frac{
\left|\sin\left(\pi \frac{\sqrt{q}}{L} \left(f_{|S|-1} + \sum
_{j=1}^{|S|-2}\epsilon_j f_j\right) \right)\right| }{
\left(f_{|S|-1} + \sum_{j=1}^{|S|-2} \epsilon_j
f_j\right)^{3/2-\delta}}
\end{multline*}
Changing sums into integrals gives
\begin{multline*}
Q(q) \ll  \idotsint_1^\infty \;\d x_1\dots\d x_{|S|-2}
\sum_{\epsilon_j = \pm 1} \int_{1+\max\{0,-\sum_{j=1}^{|S|-2}
\epsilon_j x_j\}}^{\infty} \;\d x_{|S|-1}\\
\times \left(\prod_{j=1}^{|S|-1} \frac{|\sin(\pi \tfrac{1}{L} x_j
\sqrt{q} )| }{x_j^{3/2-\delta}}\right) \frac{|\sin(\pi
\tfrac{1}{L} \sqrt{q} (x_{|S|-1}+\sum_{j=1}^{|S|-2}\epsilon_j
x_j))|}{\left(x_{|S|-1}+\sum_{j=1}^{|S|-2} \epsilon_j
x_j\}\right)^{3/2-\delta}}
\end{multline*}
and changing variables to $x_j \sqrt{q}/L \to y_j$,
\begin{multline*}
Q(q) \ll \frac{q^{|S|/4+1/2-|S|\delta/2}}{L^{|S|/2+1-|S|\delta}}
\idotsint_{\sqrt{q}/L}^\infty  \sum_{\epsilon_j = \pm 1}
\int_{\max\{0,-\sum \epsilon_j y_j\}+\sqrt{q}/L}^\infty \\
\left(\prod_{j=1}^{|S|-1} \frac{|\sin(\pi y_j )|
}{y_j^{3/2-\delta}}\right) \frac{|\sin(\pi
(y_{|S|-1}+\sum_{j=1}^{|S|-2}\epsilon_j
y_j)|}{\left(y_{|S|-1}+\sum_{j=1}^{|S|-2} \epsilon_j
y_j\right)^{3/2-\delta}} \;\d y_{|S|-1}\d y_{|S|-2} \dots\d y_{1}
\end{multline*}

Since the multiple integral is bounded, we may conclude that
\begin{equation*}
Q(q) \ll
\begin{cases}
\frac{q^{|S|/4+1/2-|S|\delta/2}}{L^{|S|/2+1-|S|\delta}} & \text{
if }
q<L^2\\
1 & \text{ if } q \geq L^2
\end{cases}
\end{equation*}
substituting this into \eqref{eq:upper bound D1} we see that
\begin{equation*}
\sum_{q \sqfree} D_q(S) \ll
\begin{cases}
\frac{L^{\delta'}}{L} & \text{ if } |S|=2\\
\frac{L^{\delta'}}{L^{|S|/2+1}} & \text{ if } |S|\geq 3
\end{cases}
\end{equation*}

Hence
\begin{equation*}
\frac{1}{\sigma^{|S|}}\sum_{q \sqfree} D_q(S) \ll
\begin{cases}
L^{\delta'} & \text{ if } |S|=2\\
\frac{1}{L^{1-\delta'}} & \text{ if } |S|\geq 3
\end{cases}
\end{equation*}
since equation \eqref{eq:sigma_asympt} gives $\sigma \sim
\frac{4\sqrt{\log L}}{\sqrt{L}}$ when $L\to\infty$ but
$L/\sqrt{M}\to 0$. However, in the case $|S|=2$, by the definition
of $D_q(S)$ and $\sigma^2$ we see that
\begin{equation*}
\sum_{q \sqfree} D_q(S) = \sigma^2.
\end{equation*}
This completes the proof of the lemma.
\end{proof}

\section{Unsmoothing}\label{sec: unsmoothing}
Recall that $S(t,1/L)$ is the normalized remainder term for the
number of lattice points in an annulus of inner radius $t$ and
width $1/L$. In this section we prove Theorem~\ref{thm:Bleher
conjecture} by showing that the variance of the difference
$(S(t,1/L)-\Sf(t))/\sigma$ vanishes and then combining this with
Chebyshev's inequality to deduce a distribution theorem for
$S(t,1/L)$.

We begin with an approximation result for $N(t)$:
\begin{lem}\label{lem:Titchmarsh}
For any $a>0$, $c>1$ we have
\begin{equation*}
N(t)= \pi t^2 -\frac{\sqrt{t}}{\pi} \sum_{n\leq X}
\frac{r(n)}{n^{3/4}} \cos(2\pi t \sqrt{n} +\tfrac14\pi) +
\O(|t|^{-1/2}) + \O\left(X^a\right) +
\O(\frac{|t|^{2c-1}}{\sqrt{X}})
\end{equation*}
\end{lem}
This Lemma was already invoked by Heath-Brown in \cite{H-B}, with
the proof being an argument similar to that which derives (12.4.4)
in \cite{Tit1}.

\begin{lem}\label{lem:Chebychev}
Suppose $L\to\infty$ as $T\to\infty$ and choose $M$ so that
$L/\sqrt{M} \to 0$ as $T\to\infty$ but $M=\O(T^{2(1-\delta)})$
(for a fixed $\delta>0$). Then as $T\to\infty$,
$$
\<|S(t,1/L)-\Sf(t)|^2 \> \ll \frac{\log M}{\sqrt{M}}
$$
\end{lem}

\begin{proof}
Putting $a=\delta'$ and $c=1+\delta'/2$ for $\delta'>0$
arbitrarily small in Lemma~\ref{lem:Titchmarsh}, we have
\begin{align*}
S(t,1/L) &:= \frac{N(t+1/L) - N(t) -
\pi(2t/L+1/L^2)}{\sqrt{t}}\\
&=\frac{2}{\pi} \sum_{n\leq X} \frac{ r(n)}{ n^{3/4}}
\sin\left(\frac{\pi \sqrt{n}}{L}\right) \sin\left(2\pi
(t+\tfrac{1}{2L})\sqrt{n}+\tfrac14\pi\right) +R(X,t)
\end{align*}
where
\begin{equation*}
R(X,t) \ll \frac{1}{|t|} + \frac{X^{\delta'}}{\sqrt{|t|}} +
\frac{|t|^{1/2+\delta'}}{\sqrt{X}}
\end{equation*}
Set $X=T^{2-\delta}$. Since $M=\O(T^{2(1-\delta)})$ and $\^\psi$
has compact support, the infinite sum in $\Sf(t)$, given in
\eqref{eq:smooth S infinite sum}, is truncated before
$n=T^{2-\delta}$, and so
\begin{multline*}
S(t,1/L) - \Sf(t) = \\
\frac{2}{\pi} \sum_{n\leq T^{2-\delta}} \frac{ r(n)}{ n^{3/4}}
\sin\left(\frac{\pi \sqrt{n}}{L}\right) \sin\left(2\pi
(t+\tfrac{1}{2L})\sqrt{n}+\tfrac \pi 4\right)
\left(1-\^\psi(\sqrt{\frac{n}{M}})\right) +R(T^{2-\delta},t)
\end{multline*}
Let $P$ denote the sum, then Cauchy-Schwartz gives
\begin{equation}\label{eq:defn P}
\<(S(t,1/L) - \Sf(t))^2\> = \<P^2\> +\<R(T^{2-\delta},t)^2\> +
\O\left(\sqrt{\<P^2\>}\sqrt{\<R(T^{2-\delta},t)^2\>}\right)
\end{equation}

Observe that
\begin{equation*}
\<R(T^{2-\delta},t)^2\> \ll T^{-1+\delta''}
\end{equation*}
for arbitrarily small $\delta''>0$, and
\begin{multline*}
\<P^2\> = \frac{2}{\pi^2} \sum_{n\leq T^{2-\delta}} \frac{
r(n)^2}{ n^{3/2}} \sin^2\left(\frac{\pi \sqrt{n}}{L}\right)
\left(1-\^\psi(\sqrt{\frac{n}{M}})\right)^2\\
+\O\left(\sum_{1\leq m\neq n\leq T^{2-\delta}}
\^\omega(T(\sqrt{n}-\sqrt{m}))\right)
\end{multline*}
The same argument used in \S\ref{sect:variance} shows the error
term here vanishes like $\O(T^{-B})$ for any $B>0$.

Since $\sum_{n\leq X} r(n)^2 \sim4X\log X$, partial summation
gives
\begin{align*}
\<P^2\> &\sim \frac{8}{\pi^2} \int_1^{T^{2-\delta}} \frac{ \log
x}{ x^{3/2}} \sin^2\left(\frac{\pi \sqrt{x}}{L}\right)
\left(1-\^\psi(\sqrt{\frac{x}{M}})\right)^2 \;\d x\\
&= \frac{32}{L \pi^2} \int_{1/L}^{T/L} \frac{\log(Ly) \sin^2(\pi
y)}{y^2} \left(1-\^\psi(\frac{y L}{\sqrt{M}})\right)^2 \;\d y
\end{align*}
by change of variables $x=y^2 L^2$. If $yL / \sqrt{M} \ll 1$ then
\begin{equation*}
\^\psi(\frac{y L}{\sqrt{M}}) = 1 + \O\left(
\frac{yL}{\sqrt{M}}\right)
\end{equation*}
leading to
\begin{align*}
\<P^2\> &\ll \frac{L}{M} \int_0^{\sqrt{M}/L} \log(Ly) \sin^2(\pi
y) \;\d y + \frac{1}{L} \int_{\sqrt{M}/L}^{T/L} \frac{ \log(Ly)
\sin^2(\pi
y)}{y^2} \;\d y\\
&\ll \frac{\log M}{\sqrt{M}}
\end{align*}
Inserting this into \eqref{eq:defn P}, using
$M=\O(T^{2(1-\delta)})$ and choosing $0<\delta''<\delta$ in the
estimate of $\<R(X,t)^2\>$ we have that
\begin{align*}
\<(S(t,1/L) - \Sf(t))^2\> &\ll\frac{\log M}{\sqrt{M}} +
\frac{1}{T^{1-\delta''}} + \frac{\sqrt{\log M}}{M^{1/4}
T^{1/2-\delta''/2}}\\
&= \O\left(\frac{\log M}{\sqrt{M}} \right)
\end{align*}
\end{proof}

\begin{lem}\label{lem: close ditributions}
Under the conditions of Lemma~\ref{lem:Chebychev}, we have for all
fixed $\eta>0$,
\begin{equation*}
\Prob_{\omega,T}\left\{ \left|\frac{S(t,1/L)}{\sigma} -
\frac{\Sf(t)}{\sigma}\right|
> \eta \right\}  \to 0
\end{equation*}
as $T\to\infty$, where $\sigma^2 = 16\log L / L$.
\end{lem}
\begin{proof}
For fixed $\eta>0$, Chebychev's inequality gives
\begin{align*}
\Prob_{\omega,T}\left\{ \left|\frac{S(t,1/L)}{\sigma} -
\frac{\Sf(t)}{\sigma}\right|
> \eta \right\} &\leq \frac{\<(S(t,1/L) - \Sf(t))^2\>}{\eta^2 \sigma^2}\\
&\ll \frac{L}{\log L} \frac{\log M}{\sqrt{M}}
\end{align*}
which tends to zero as $T\to\infty$ by the assumptions placed on
$M$ and $L$.
\end{proof}

\begin{cor}\label{lem:S0 Gaussian, smooth average}
If $L\to\infty$ but $L=\O(T^\delta)$ for all $\delta>0$ as
$T\to\infty$, then for any interval $\mathcal{A}$,
\begin{equation}
\Prob_{\omega,T}\left\{\frac{S(t,1/L)}{\sigma} \in
\mathcal{A}\right\} \to \frac{1}{\sqrt{2\pi}} \int_{\mathcal{A}}
e^{-x^2/2}\;\d x
\end{equation}
where $\sigma^2 = 16\log L / L$.
\end{cor}

\begin{proof}
Set $M=L^3$, then $M=\O(T^\delta)$ for all $\delta>0$ and
$L/\sqrt{M} \to 0$. Thus $\Sf/\sigma$ weakly converges to a
standard normal distribution as $T\to\infty$ when $t$ is smoothly
averaged around $T$ by Theorem~\ref{thm:distrn for Sf}. But
Lemma~\ref{lem: close ditributions} implies that $S(t,1/L)/\sigma$
must also weakly converge to a standard normal distribution too.
\end{proof}

We are now able to prove our main result, Theorem~\ref{thm:Bleher
conjecture}, which says that if $L \to \infty$ but $L=
\O(T^\delta)$ for all $\delta>0$, then for any  interval
$\mathcal{A}$
\begin{equation*}
\lim_{T\to\infty} \frac{1}{T}\meas\left\{t\in[T,2T]\ :\
\frac{S(t,1/L)}{\sigma} \in \mathcal{A}\right\} =
\frac{1}{\sqrt{2\pi}} \int_\mathcal{A} e^{-x^2/2}\;\d x
\end{equation*}

\begin{proof}[Proof of Theorem~\ref{thm:Bleher conjecture}]
Fix $\epsilon>0$, and approximate the indicator function
$\I_{[1,2]}$ above and below by smooth functions $\chi_{\pm}\geq
0$ so that $\chi_{-}\leq \I_{[1,2]} \leq \chi_{+}$, where both
$\chi_{\pm}$ and their Fourier transforms are smooth  and of rapid
decay, and so that their total masses are within $\epsilon$ of
unity: $| \int\chi_{\pm}(x)\d x-1 |<\epsilon$. Now set
$\omega_{\pm}:=\chi_{\pm}/\int \chi_{\pm}$. Then $\omega_{\pm}$
are  ``admissible''  and for all $t$,
\begin{equation}\label{approx}
(1-\epsilon)\omega_{-}(t) \leq \I_{[1,2]}(t) \leq (1+\epsilon)
\omega_{+}(t)
\end{equation}

Now
$$
\meas\left\{t\in[T,2T]\ :\ \frac{S(t,1/L)}{\sigma} \in
\mathcal{A}\right\} =\int_{-\infty}^\infty
\I_{\mathcal{A}}\left(\frac{S(t,1/L)}{\sigma}\right)
\I_{[1,2]}\left(\frac tT\right) dt
$$
and since \eqref{approx} holds, we find
\begin{multline*}
(1-\epsilon)
\Prob_{\omega_{-},T}\left\{\frac{S(t,1/L)}{\sigma}\in\mathcal{A}\right\}
\leq \frac{1}{T} \meas\left\{t\in[T,2T]\ :\
\frac{S(t,1/L)}{\sigma} \in
\mathcal{A}\right\}\\
\leq (1+\epsilon) \Prob_{\omega_{+},T}
\left\{\frac{S(t,1/L)}{\sigma}\in\mathcal{A}\right\}
\end{multline*}
By Corollary~\ref{lem:S0 Gaussian, smooth average} the two extreme
sides of this inequality have a limit as $T\to \infty$, of
$$ (1 \pm \epsilon) \frac{1}{\sqrt{2\pi}}\int_{\mathcal{A}} e^{-x^2/2} \;\d x $$
and so we get that
$$
(1-\epsilon)\frac{1}{\sqrt{2\pi}}\int_{\mathcal{A}} e^{-x^2/2}
\;\d x  \leq \liminf_{T\to\infty} \frac{1}{T}
\meas\left\{t\in[T,2T]\ :\ \frac{S(t,1/L)}{\sigma} \in
\mathcal{A}\right\}
$$ with a similar statement for limsup; since $\epsilon>0$ is
arbitrary this shows that the limit exists and equals
$$\lim_{T\to\infty} \frac 1T \meas\left\{t\in[T,2T]\ :\ \frac{S(t,1/L)}{\sigma}
\in
\mathcal{A}\right\}  = \frac{1}{\sqrt{2\pi}}\int_{\mathcal{A}}
e^{-x^2/2} \;\d x$$ which is the Gaussian law.
\end{proof}

\section*{Acknowledgments}

We would like to thank Jens Marklof for helpful discussions. This
work was supported in part by the EC TMR network
\textit{Mathematical Aspects of Quantum Chaos}, EC-contract no
HPRN-CT-2000-00103 and the Israel Science Foundation founded by
the Israel Academy of Sciences and Humanities.

\end{document}